\documentclass{amsart}

\usepackage{amsfonts}
\usepackage{cases}
\usepackage{mathrsfs}
\usepackage{bbm}
\usepackage{amssymb}
\usepackage{txfonts}
\usepackage{amscd}
\usepackage{amsfonts,latexsym,amsmath,amsthm,amsxtra,mathdots}
\usepackage[bookmarks,colorlinks]{hyperref}
\usepackage[all,cmtip]{xy}
\RequirePackage{amsmath} \RequirePackage{amssymb}
\usepackage{color}
\usepackage{colordvi}
\usepackage{multicol}
\usepackage{tikz}
\usepackage{hyperref}
\usepackage{graphicx}
\usepackage{amsmath}
\usepackage{amsmath,amscd}
\usepackage[normalem]{ulem}
\usetikzlibrary{matrix,arrows}
\usepackage{textcomp}
\usepackage{mathtools}

\usepackage{cleveref}

\newcommand{\id}{{{\rm id}}}

\newcommand{\subgroup}{\leq}

    \newcommand{\BE}{{\mathbb {E}}}

     \newcommand{\BN}{{\mathbb {N}}}
     
    \newcommand{\BQ}{{\mathbb {Q}}} \newcommand{\BR}{{\mathbb {R}}}

     \newcommand{\BZ}{{\mathbb {Z}}}

    \newcommand{\im}{{\mathrm{im}}} 
     \newcommand{\rank}{{\mathrm{rank}}}

\def\-{^{-1}}

\newcommand{\delete}[1]{}

    \theoremstyle{plain}

\newtheorem{thm}{Theorem}[section]
\newtheorem{defn}[thm]{Definition} 
\newtheorem{defns}[thm]{Definitions}
\newtheorem{ex}[thm]{Example} 
\newtheorem{lem}[thm]{Lemma}
\newtheorem{prop}[thm]{Proposition}
\newtheorem{cor}[thm]{Corollary}
\newtheorem{rem}[thm]{Remark}

\newtheorem*{thm*}{Theorem}

\newtheorem*{rem*}{Remark}
\newtheorem*{conjA}{Conjecture A}
\newtheorem*{conjB}{Conjecture B}

    \numberwithin{equation}{section}

\newcommand{\uueg}{{\underline{\underline E}}G}

\newcommand{\uu}[1]{\uline{\uline{#1}}}
\renewcommand{\u}[1]{\uline{#1}}

\def\Proof{\noindent{\bf Proof}\quad}
\def\qed{\hfill$\square$\smallskip}

\newenvironment{cond}[1]{%
  \vspace{1.2ex}\begin{enumerate}%
  \renewcommand*{\theenumi}{#1}%
  \renewcommand*{\labelenumi}{#1}%
  \item%
}{%
  \end{enumerate}\vspace{1.2ex}%
}

\begin{document}

\title{On the finiteness of the classifying space for the family of virtually cyclic subgroups}

\author{Timm von Puttkamer}
\address{University of Bonn, Mathematical Institute, Endenicher Allee 60, 53115 Bonn, Germany}
\email{tvp@math.uni-bonn.de}

\author{Xiaolei Wu}
\address{Max Planck Institut f\"ur Mathematik, Vivatsgasse 7, 53111 Bonn, Germany }
\email{hsiaolei.wu@mpim-bonn.mpg.de}


\date{July, 2017}


\begin{abstract}
Given a group $G$, we consider its classifying space for the family of virtually cyclic subgroups. We show for many groups, including for example, one-relator groups, acylindrically hyperbolic groups, $3$-manifold groups  and CAT(0) cube groups, that they do not admit a finite model for this classifying space unless they are virtually cyclic. This settles a conjecture due to Juan-Pineda and Leary for these classes of groups.
\end{abstract}

\maketitle
\section*{Introduction}

Given a group $G$, we denote by $\uueg = E_{\mathcal{VCY}}(G)$ a $G$-CW-model for the classifying space for the family of virtually cyclic subgroups. The space $\uueg$ is characterized by the property that the fixed point set $\uueg^H$ is non-empty and contractible for all virtually cyclic subgroups $H$ of $G$ and empty otherwise. In the formulation of the Farrell--Jones Conjecture $\uueg$ plays an important role (see for example \cite{FJ,LuRe} for more information). Due to this, there has been a growing interest in studying $\uueg$, see for example \cite{DP, JL, KMN,Lu05, Lu09}.  Recall that a $G$-CW-complex $X$ is said to be finite if it has finitely many orbits of cells. Similarly, $X$ is said to be of finite type if it has finitely many orbits of cells of dimension $n$ for any $n$. In \cite[Conjecture 1]{JL}, Juan-Pineda and Leary formulated the following conjecture:

\begin{conjA}   \cite [Juan-Pineda and Leary]{JL}
Let $G$ be a group admitting a finite model for $\uueg$. Then $G$ is virtually cyclic.
\end{conjA}

The conjecture may be surprising at the beginning as there are many examples of groups with a finite model for the classifying space for the family consisting only of the trivial subgroup or for the family of finite subgroups. Juan-Pineda and Leary demonstrated the validity of the conjecture for hyperbolic groups \cite[Corollary 12]{JL}. Later, Kochloukova, Mart\'inez-P\'erez and Nucinkis verified the conjecture for elementary amenable groups \cite{KMN}.  Groves and Wilson gave a simplified proof for elementary amenable groups \cite{GW}. So far, almost all the proofs boil down to analyzing whether $\uueg$ has a finite $0$-skeleton. It turns out that having a finite $0$-skeleton for $\uueg$ is equivalent to the following purely algebraic condition (see \Cref{BVCfinitezeroskeleton})

\begin{cond}{(BVC)}\label{cond:BVC}
\textit{ $G$ has a finite set of virtually cyclic subgroups $\{V_1,V_2,\ldots, V_n\}$ such that every virtually cyclic subgroup of $G$ is conjugate to a subgroup of some $V_i$.}
\end{cond}

Following Groves and Wilson, we shall call this property BVC and the finite set $\{V_1,V_2,\\\ldots, V_n\}$ a \textbf{witness} to BVC for $G$. In this paper, we give a systematic study of the property BVC. Our main theorem can be stated as follows

\begin{thm*}
The following classes of groups satisfy Conjecture A:
\begin{enumerate}

    \item HNN extensions of free groups of finite rank (\ref{HNN-ext-of-free-groups-BVC}),
    \item one-relator groups (\ref{one-relator-BVC}),
    \item acylindrically hyperbolic groups (\ref{acylindrically hyperbolic}),
    \item $3$-manifold groups (\ref{3-manifold-group-BVC}),
    \item CAT(0) groups that contain a rank-one isometry or $\BZ^2$ (\ref{BVCCAT(0)}), in particular CAT(0) cube groups (\ref{BVCCAT(0)cube}).
\end{enumerate}

\end{thm*}

In fact, we show that any finitely generated group in these classes has BVC if and only if it is virtually cyclic.  Our result also suggests that the following conjecture could be true.

\begin{conjB}\label{conj-bvc}
Let $G$ be a finitely presented group which has BVC. Then $G$ is virtually cyclic.
\end{conjB}

\begin{rem*}

\begin{enumerate}
\item The assumption of having a finitely presented group is necessary here since
Osin \cite{Os10} has constructed finitely generated torsion-free groups with exactly two conjugacy classes. In particular these groups have BVC.

\item We do not know whether Conjecture B holds for all elementary amenable groups. Groves and Wilson showed that solvable groups satisfy it \cite{GW}.
\item If we knew the Flat Closing Conjecture, then it would follow that any CAT(0) group satisfies Conjecture B. See \Cref{flat-closing-conj} for more information.
\end{enumerate}
\end{rem*}

Our paper is organized as follows. In \Cref{basic}, we first summarize what we already know about groups admitting a finite model for $\uueg$, then we study basic properties of groups with BVC and deduce that many groups cannot have this property. In \cref{HNNonerelator}, we study HNN extension of groups and show for example that any HNN extension of a finitely generated free group does not have BVC. Using this, we prove that any non-cyclic one-relator group does not have BVC. In \cref{hyper}, we show that acylindrically hyperbolic groups and finitely generated non virtually cyclic $3$-manifold groups do not have BVC.  In the last section, we study groups acting on CAT(0) spaces. We show for example, that CAT(0) cube groups do not have BVC unless they are virtually cyclic.

\vspace{3mm}

Results of this paper will also appear as a part of the first author's thesis.

\textbf{Acknowledgements.} The first author was supported by an IMPRS scholarship of the Max Planck Society. The second author would like to thank the Max Planck Institute for Mathematics at Bonn for its support. We also want to thank Yongle Jiang for helpful comments.

\section{Properties of groups admitting a finite model for $\uueg$}\label{basic}
In this section we first review properties of groups admitting a finite model for $\uueg$. We then proceed to prove many useful lemmas for groups with BVC. We also use these lemmas to show many groups cannot have BVC.

We denote by $\underline{E}G$, resp.~$EG$ the $G$-CW classifying space for the family of finite subgroups resp.~for the family consisting only of the trivial subgroup. We summarize the properties of groups admitting a finite model for $\uueg$ as follows

\begin{prop}
Let $G$ be a group admitting a finite model for $\uueg$, then

\begin{enumerate}
    \item $G$ has BVC.

    \item $G$ admits a finite model for $\underline{E}G$.

    \item For every finite subgroup of $H \subset G$, the Weyl group $W_GH$ is finitely presented and of type $FP_{\infty}$. Here $W_GH := N_G(H)/H$, where $N_G(H)$ is the normalizer of $H$ in $G$.

    \item $G$ admits a finite type  model for $EG$. In particular, $G$ is finitely presented.

    \item $G$ has finitely many conjugacy classes of finite subgroups. In particular, the order of finite subgroups of $G$ is bounded.

\end{enumerate}

\end{prop}

\Proof Note that if $G$ has a finite or finite type  model for $\underline{E}G$, then $G$ is finitely presented and has finitely many conjugacy classes of finite subgroups \cite[4.2]{Lu00}. Now (a), (b), (c) can be found for example in \cite[Section 2]{KMN}. Part (d) can be deduced from (c) by taking $H$ to be the trivial group.

\begin{rem}
If one replaces finite by finite type in the assumptions of the above proposition, then the conclusions still hold except one has to replace finite by finite type in (b).
\end{rem}

The following lemma is well-known to experts.

\begin{lem} \label{BVCfinitezeroskeleton}
Let $G$ be a group. Then there is a model for $\uueg$ with finite $0$-skeleton if and only if $G$ has BVC.
\end{lem}

\Proof Suppose $X$ is a $G$-CW model for $\uueg$ with a finite $0$-skeleton. Let $\sigma_1,\sigma_2,\ldots, \sigma_n$ be a set of representatives for each orbit of $0$-cell. Let $V_1,V_2,\ldots,V_n$ be the corresponding virtually cyclic stabilizers. Since for every virtually cyclic subgroup $V$, the set of fixed points $X^V$ is non-empty, there exists some vertex of $X$ that is fixed by $V$. Since this vertex stabilizer is a conjugate to some $V_i$, the subgroup $V$ is subconjugate to $V_i$. Conversely, suppose $G$ has BVC, and let $V_1,V_2,\ldots,V_n$ be witnesses. We can construct a model for $\uueg$ with finite $0$-skeleton as follows. Consider the $G$-set $S:= \coprod_{i = 1}^n G/V_i$. The complete graph $\Delta(S)$ spanned by $S$ , i.e., the simplicial graph that
contains an edge for every two elements in  $\Delta(S)$, carries a cocompact simplicial $G$-action.  The first barycentric subdivision  $\bar{\Delta}(S)$ of  $\Delta(S)$  is a $G$-CW-complex. Note that $\bar{\Delta}(S)$ is again cocompact. Moreover, $\bar{\Delta}(S)^H$ is nonempty when $H$ is virtually cyclic and empty otherwise. Now we can add equivariant cells of dimension $\geq 1$ to $\bar{\Delta}(S)$ and make $\bar{\Delta}(S)^H$ contractible for all virtually cyclic subgroup by induction using \cite[Proposition 2.3]{Lu89}. This way we obtain a model for $\uueg$ with finite $0$-skeleton.
\qed

The following structure theorem about virtually cyclic groups is well known, see for example \cite[Proposition 4]{JL} for a proof.

\begin{lem}  \label{vcstru}
Let $G$ be a virtually cyclic group. Then $G$ contains a unique maximal normal finite subgroup $F$ such that one of the following holds

\begin{enumerate}

\item the finite case, $G =F$;
\item the orientable case, $G/F$ is the infinite cyclic group;
\item  the nonorientable case, $G/F$ is the infinite dihedral group.

\end{enumerate}

\end{lem}

Note that the above lemma implies that a torsion-free virtually cyclic group is either trivial or isomorphic to $\BZ$. Thus we have the following
\begin{cor}\label{torsionfreeBVC}
Let $G$ be a torsion-free group, then $G$ has BVC if and only if there exist elements $g_1,g_2,\ldots g_n$ in $G$ such that every element in $G$ is conjugate to a power of some $g_i$.
\end{cor}

The following lemma is key to many of our proofs.

\begin{lem} \label{vircycorder}
Let $V$ be a virtually cyclic group and let $g,h \in V$ be two elements of infinite order, then there exist $p,q \in \BZ$ such that $g^p =h^q$. Furthermore, there exists $v_0 \in V$ such that for any $v\in V$ of infinite order there exist nonzero $p_0,p$ such that  $v_0^{p_0} = v^p$ with $\frac{p_0}{p} \in \BZ$.
\end{lem}

\Proof By \Cref{vcstru}, there exists an finite normal subgroup of $F$ such that $V/F$ is isomorphic to $\BZ$ or $\BZ \rtimes \BZ/2$. We denote the quotient map by $q$. Then if $g$ is of infinite order and $V/F \cong \BZ \rtimes \BZ/2$, then $q(g) \in \BZ \rtimes \{0\}$. Thus we can always assume $V/F \cong \BZ$. In this case $V \cong F \rtimes_f \BZ$, where $\BZ$ acts on $F$ via the automorphism $f$. Now for any $g = (x,m) \in F \rtimes_f \BZ$ we can choose $k = |F||f|$ where $|F|$ is the order of the finite group $F$ and $|f|$ is the order of the automorphism $f$, then $g^k = (0,km)$. Let $h=(y,n)$, then $h^k = (0,kn)$. Now we see that $g^{kn} =h^{km}$. If we choose $v_0 = (0,1) \in F\rtimes_f \BZ$ then for any $v= (x,m) \in V$,  $v_0^{mk} = v^{k}$.
\qed

Note that in any virtually cyclic group there are only finitely many distinct finite subgroups up to conjugacy (using \Cref{vcstru} or the fact that virtually cyclic groups are CAT(0)). Using this fact one immediately obtains

\begin{lem} \label{BVCfinitesubgroup}
If a group $G$ has BVC, then $G$ has finitely many conjugacy classes of finite subgroups. In particular, the order of finite subgroups in $G$ is bounded.
\end{lem}

In a group $G$, we call an element $g$ \textbf{primitive} if it cannot be written as a proper power. Then \Cref{torsionfreeBVC} implies the following

\begin{lem} \label{primitiveBVC}
Let $G$ be a torsion-free group. If $G$ has infinitely many conjugacy classes of primitive elements, then $G$ does not have BVC.
\end{lem}

Note that without the assumption of $G$ being torsion-free, the previous lemma does not hold. In fact, even a virtually cyclic group could contain infinitely many conjugacy classes of primitive elements.

\begin{ex}
Let $\text{S}_n$ be the symmetric group of order $n$ with $n \geq 2$, then $\text{S}_n \times \BZ$ has infinitely many primitive conjugacy classes. In fact, let $(x_i,2^i)\in S_n \times \BZ$ with $x_i$ an odd element in $\text{S}_n$. Then $(x_i,2^i)$ is primitive for any $i \geq 1$. Since if $(x_i,2^i) = (x,y)^k$ for some $k>1$, then $ky=2^i$. In particular, $k$ is even and thus $x^k$ cannot equal the odd element $x_i$. 
The elements $(x_i, 2^i)$ cannot be conjugate to each other since the second coordinate is different.
\end{ex}

\begin{lem} \label{free-products-bvc}
Let $G = A\ast B$ be a free product with $A$ and $B$ nontrivial, then $G$ has BVC if and only if $G$ is virtually cyclic.
\end{lem}
\Proof If $A$ and $B$ are finite groups, then $A \ast B$ is a virtually free group and hence hyperbolic. So the lemma holds in this case. Now we assume that $A$ is not finite, then $A \ast B$ is not virtually cyclic. Let $a_1,a_2 \ldots, a_n,\ldots,$ be a sequence of different elements in $A$ and let $b \in B$ be a non-trivial element. Then $\{a_ib \mid i\geq1\}$ in $G$ form infinitely many conjugacy classes of primitive elements in $G$. Moreover, when $i \neq j$, $a_ib$ and any conjugates of $a_jb$ cannot lie in a virtually cyclic group. In fact, if this were the case, by \Cref{vircycorder} we would have that $(a_ib)^m$ is conjugate to $(a_jb)^n$ in $G$ for some $m,n \neq 0$. But this is impossible by the choices of $a_i$ and $b$ and \cite[IV.2.8]{LS77}. Hence $G$ does not have BVC in this case.
\qed

\begin{lem} \cite[5.6]{KMN} \label{finiteindexsubgroupbvc}
If a group $G$ has BVC, then any finite index subgroup also has BVC.
\end{lem}

Combining this with the main result of \cite{GW}, we have

\begin{prop} \label{BVCvirtuallysolvable}
Virtually solvable groups have BVC if and only if they are virtually cyclic.
\end{prop}

\begin{lem}\cite[2.2]{GW} Let $G$ be a group with property BVC. Then the following assertions hold.

\begin{enumerate}\label{bvc-properties}
    \item The group $G$ satisfies the ascending chain condition for normal subgroups.
    \item If $L$ and $M$ are normal subgroups of $G$ with $M < L$ and $L/M$ a torsion group, then there are only finitely many normal subgroups  $K$ of $G$ such that $M \leq K \leq L$. \label{bvc-properties-quotients}
    \item Let
    $$1 = G_n \leq G_{n-1} \leq \cdots \leq G_1 \leq G_0 = G$$
    be a series of normal subgroups of $G$. Then the number of factors $G_i/G_{i-1}$ that are not torsion groups is bounded by the number of infinite groups in a witness to $BVC$ for $G$.
\end{enumerate}
\end{lem}

The following lemma allows us to show that many groups cannot have BVC.

\begin{lem}\label{BVCquotient}
Let $\pi \colon G \rightarrow Q$ be a surjective group homomorphism. If $Q$ is a torsion-free group that does not have BVC, then $G$ does not have BVC.
\end{lem}

\Proof Suppose $G$ has BVC and let $V_1,V_2\ldots, V_n$ be a witness for BVC of $G$. Note that any quotient of a virtually cyclic group is again virtually cyclic. Hence $\pi(V_i)$ are virtually cyclic subgroups in $Q$. Since $Q$ is torsion-free and does not have BVC, we can find a nontrivial element $q \in Q$ such that $q$ cannot be conjugated to $\pi(V_i)$ for any $i$. Now take $g \in G$ such that $\pi(g) = q$, then we can find $c \in G$ such that $cgc^{-1} \in V_i$ for some $i$. But then we would have $\pi(c)\pi(g)\pi(c^{-1}) =  \pi(c) q \pi(c^{-1})\in \pi(V_i)$ which is a contradiction.    \qed

\begin{cor}\label{BVChomology}
If $G$ is a group having BVC, then the abelianization $H_1(G, \BZ)$ is finitely generated of rank at most one.
\end{cor}
\Proof
Let $A$ be the abelianization of $G$ and let $T$ be the torsion subgroup of $A$. By \Cref{BVCquotient} the torsion-free abelian group $A/T$ has to be trivial or infinite cyclic. Hence $A = T$ or $A \cong T \times \BZ$. Now \Cref{bvc-properties}(b) implies that the torsion group $T$ has to be finite.
\qed

\begin{ex}
The Thompson groups are a family of three finitely presented groups, $F \subset T\subset V $. Thompson's group $F$ can be defined by the presentation $\langle A,B \mid [AB^{-1},A^{-1}BA] = [AB^{-1},A^{-2}BA^2]= 1 \rangle$. Since $H_1(F) \cong \BZ^2$, it follows that $F$ does not have BVC. Since the order of finite subgroups in $T$ and $V$ is unbounded, we see that $T$ and $V$ also do not have BVC. See \cite{CFP} for more information about Thompson groups.
\end{ex}

Recall that a group satisfies the \textbf{strong Tits alternative} if any finitely generated subgroup has a finite index subgroup which is either  solvable, or has a non-abelian free quotient. Since virtually solvable groups and free groups have BVC if and only if they are virtually cyclic, we have the following by \Cref{finiteindexsubgroupbvc} and \Cref{BVCquotient},

\begin{lem}
If a group $G$ satisfies the strong Tits alternative, then a finitely generated subgroup of $G$ has BVC if and only if it is virtually cyclic.
\end{lem}
~\qed

Since Coxeter group and right angled Artin groups are known to satisfy the strong Tits alternative by \cite{NV} and \cite[1.6]{AM}, we have the following

\begin{cor}
Let $G$ be a Coxeter group or a right angled Artin group, then a finitely generated subgroup of $G$ satisfies BVC if and only if it is virtually cyclic.
\end{cor}

Note that \Cref{BVChomology} reduces the study of finitely generated groups with BVC to the case where $H_1(G,\BZ)$ is of rank one or zero. If $H_1(G,\BZ)$ is of rank one, then $G$ surjects to $\BZ$ and thus $G$ becomes a semidirect product of the form $H \rtimes \BZ$. We proceed to study groups of this type.

Given an automorphism $\phi$ of $G$, we say that two elements $g,h$ in $G$ are $\phi$-conjugate if $g = x h \phi({x^{-1}})$ for some $x\in G$. This is an equivalence relation whose equivalence classes are called $\phi$-twisted conjugacy classes. For $\phi = \id_G$ one recovers the usual notion of conjugacy.

\begin{lem}\label{twistBVC}
Let $\phi$ be an automorphism of $H$ such that $H$ has infinitely many $\phi$-twisted conjugacy classes, then the semidirect product $ H \rtimes_\phi \BZ$ does not have BVC.
\end{lem}

\Proof Note that in $H \rtimes_\phi \BZ$, the elements $(g,1)$ and $(h,1)$ are primitive and they are in the same conjugacy class if and only if $g$ and $h$ are in the same $\phi$-twisted conjugacy class in $H$.  In fact, $(g,1)$ is conjugate to $(h,1)$ in $ H \rtimes_\phi \BZ$ if and only if we can find $(x,k)$ such that $(x,k)(g,1)(x,k)^{-1} = (x\phi^k(g)\phi(x^{-1}),1)=(h,1)$. This is equivalent to saying that $g$ is $\phi$-conjugate to $\phi^k(g)$ in $H$. But $g$ and $\phi(g)$ are $\phi$-conjugate in $H$ since  $\phi(g) = g^{-1}g\phi(g)$. Hence $g$ and $h$ are in the same $\phi$-twisted conjugacy class in $H$ if and only if $(g,1)$ is conjugate to $(h,1)$ in $H \rtimes_\phi \BZ$.

Since $H$ has infinitely many $\phi$-twisted conjugacy classes, we have infinitely many primitive elements of the form $(g,1) \in H \rtimes_\phi \BZ$ that are not conjugate to each other. If $H\rtimes_\phi \BZ$ has BVC, then we can choose infinitely many elements $(g_1,1),(g_2,1),\ldots,(g_n,1),\ldots$ that are not conjugate to each other, but they lie in the same virtually cyclic subgroup. In particular the group $V$ generated by $(g_1,1),(g_2,1),\ldots,(g_n,1)\ldots$ in $G$ is virtually cyclic. But this is impossible. Consider the quotient map $q \colon H \rtimes \BZ_\phi \rightarrow \BZ$, which is onto when restricted to $V$. Hence the kernel must be finite. This means $V \cong F \rtimes \BZ$ for some finite group $F$ and the image of $(g_i,1)$ is of the form $(f_i,1)$. This leads to a contradiction since there are infinitely many $(g_i,1)$ but only finitely many $(f_i,1)$.\qed

Recall that a group is said to have property \textbf{$R_{\infty}$} if it has infinitely many $\phi$-twisted conjugacy classes for any automorphism $\phi$.

\begin{cor}\label{semi-product}
Let $G$ be a group with the property $R_\infty$, then any semidirect product $H \rtimes_\phi \BZ$ does not have BVC.
\end{cor}
Note that there are many groups with property $R_{\infty}$, for example hyperbolic groups that are not virtually cyclic, relatively hyperbolic groups and most generalzied Baumslag-Solitar groups. For more information about groups with property $R_{\infty}$ and further examples, see~\cite{FT}.

\section{HNN extension of groups and one-relator groups}\label{HNNonerelator}

In this section, we study HNN extension of groups. We first give a quick review of basic results concerning the structure of these types of groups. Then we prove that any non-ascending HNN extension of a group and any HNN extension of a finitely generated free group do not have BVC. Using this, we are able to show that one-relator groups have BVC if and only if they are virtually cyclic.

Recall that given a group $H$ and an isomorphsim $\theta \colon A \rightarrow B$ between two subgroups $A$ and $B$ of $H$, we can define a new group $H\ast_{\theta}$, called the HNN extension of $H$ along $\theta$, by the presentation $\langle H,t \mid txt^{-1} =\theta(x), x\in A \rangle $. In the study of HNN extensions of groups, Britton's Lemma and Collins' Lemma play an important rule. We give a quick review of the two lemmas and refer to \cite[IV.2]{LS77} for more details.

\begin{defn}
A sequence $g_0 , t^{\epsilon_1} , g_1 , \ldots, t^{\epsilon_n}, g_n$ of elements with $g_i \in H$ and $\epsilon_i \in \{-1,+1\}$ is said to be \textbf{reduced} if there is no consecutive sequence $t, g_i, t^{-1}$ with $g_i \in A$ or $t^{-1}, g_j,t$ with $g_j \in B$.
\end{defn}

\begin{lem} [Britton's Lemma]
If the sequence $ g_0 , t^{\epsilon_1} , g_1 , \ldots, t^{\epsilon_n} , g_n $ is reduced and $n\geq 1$, then $g_0  t^{\epsilon_1}  g_1  \cdots t^{\epsilon_n}  g_n \neq 1$ in $H\ast_{\theta}$.
\end{lem}

In the following we will not distinguish between a sequence of words as above and the element it defines in the HNN extension $H\ast_{\theta}$.

Give any $g \in H\ast_{\theta}$, we can write $g$ in a reduced form. Let $w = g_0  t^{\epsilon_1}  g_1  \cdots t^{\epsilon_n}  g_n \neq 1$ be any reduced word in  $H\ast_{\theta}$ which represents $g$. The we define the \textbf{length}  of $g$, written as $|g|$, to be the number $n$ of occurences of $t^{\pm}$ in $w$. Moreover we call an element  $w = g_0  t^{\epsilon_1}  g_1  \cdots t^{\epsilon_n}  g_n \neq 1$ cyclically reduced if all cyclic permutations of the sequence $ g_0 , t^{\epsilon_1} , g_1 , \ldots, t^{\epsilon_n} , g_n $ are reduced. Every element of $H\ast_{\theta}$ is conjugate to a cyclically reduced element.

\begin{lem} [Collins' Lemma]
Let $G = \langle H,t \mid txt^{-1} =\theta(x), x\in A \rangle$ be an HNN extension. Let $u = g_0  t^{\epsilon_1}  g_1  \cdots t^{\epsilon_n}  g_n$ and $v$ be cyclically reduced elements of $G$ that are conjugate, $n\geq 1$. Then $| u | = | v |$, and $u$ can be obtained from $v$ by taking a suitable cyclic permutation $v^{\ast}$ of $v$, which ends in $t^{\epsilon_n}$, and then conjugating by an element $z$, where $z\in A$ if $\epsilon_n = 1$, and $z \in B$ if $\epsilon_n =-1$.
\end{lem}

We are ready to prove the following:

\begin{lem} \label{HNNindexbvc}
Let $H$ be a group and let $\theta \colon A \rightarrow B$ be an isomorphism between two subgroups of $H$. If $[H:A], [H:B] \geq 2$, then the corresponding HNN extension $ G = H\ast_{\theta}$ does not have BVC.
\end{lem}

\Proof Choose $\alpha \in H \setminus B$ and $ \beta \in H \setminus A$, define
$$w_n = t^{-1}\alpha t^{n+1}\beta$$
for $n\geq 1$. Note that the elements $w_n$ are of infinite order and cyclically reduced. By Collins' Lemma, they are not conjugate to each other. If $G$ had BVC, there would exist a virtually cyclic subgroup $V \subseteq G$ such that there is an infinite subsequence $\{w_{n_i}\}$ with each $w_{n_i}$ conjugate to an infinite order element of $V$. But this cannot happen as $w_n^{p_n}$ is not conjugate to $w_m^{p_m}$ for any $p_n,p_m \neq 0$ when $n \neq m$. In fact, first note that $w_n^{p_n}$ is cyclically reduced for any $n$. So if $w_n^{p_n}$ is conjugate to $w_m^{p_m}$, their lengths must coincide by Collins' Lemma. Hence we have an equation $|p_n| (n+2)= |p_m|(m+2)$. On the other hand, there is a canonical quotient map $q \colon G \rightarrow \langle t\rangle \cong \BZ$. If $w_n^{p_n}$ is  conjugate to $w_m^{p_m}$, then $q(w_n^{p_n}) = q(w_m^{p_m})$. This means $p_n n = p_m m$. But the two equations can never hold at the same time when $n,m\geq 1$ unless $n=m$. Hence we have a contradiction by \Cref{vircycorder}.
\qed

Now when $H = A$ or $H = B$, we would have an ascending HNN extension. It seems to us that this case cannot be handled as easily as before.
In the following we will analyze the case of an ascending HNN extension of a free group $F$ of finite rank in detail. We will first deal with the case that $\theta \colon F \rightarrow F$ is injective with its image lying in the commutator subgroup of $F$. Given a group $G$, we will write $G'$ for the commutator subgroup $[G,G]$ and we will denote the $r$-th term in the lower central series by $\Gamma_r(G) = [\Gamma_{r-1}(G),G]$ where $\Gamma_1(G) = G$. Let us first recall the following facts:

\begin{lem}\label{centralseries}
The lower central series has the following properties:
\begin{enumerate}
    \item $[\Gamma_r(G),\Gamma_s(G)] \subseteq \Gamma_{r+s}(G)$ for any group $G$.
    \item $\bigcap_{r \geq 1} \Gamma_r(F) = \{ 1 \}$ for any free group $F$. 
    \item $\Gamma_{r}(F)/\Gamma_{r+1}(F)$ is a free abelian group for any $r$ and any free group $F$ of finite rank.
\end{enumerate}
\end{lem}
\Proof (a) can be found in \cite[Corollary 10.3.5]{Ha}. (b) and (c) can be found in \cite[Chapter 11]{Ha}.  \qed

\begin{cor} \label{derivedrank}
Let $\theta\colon F\rightarrow F$ be an injective map of the finitely generated free group $F$ with the image of $\theta$ lying in the commutator subgroup of $F$. If $x \in \Gamma_r(F)$, then $\theta(x) \in \Gamma_{2r}(F)$.
\end{cor}

\Proof If $x \in \Gamma_1(F) = F$, then by assumption on $\theta$ we have $\theta(x) \in [F,F] = \Gamma_2(F)$. Let $r \geq 2$ and suppose that for any $s < r$ the claim holds. If $x \in \Gamma_r(F) = [\Gamma_{r-1}(F), F]$, then by induction and \Cref{centralseries} we get $\theta(x) \in [\Gamma_{2(r-1)}(F), \Gamma_2(F)] \subgroup \Gamma_{2r}(F)$.
\qed

\begin{lem}\label{HNNextensionkernel}
Let $G = \langle H,t \mid txt^{-1} =\theta(x), x\in H \rangle$ be an ascending HNN extension of a group $H$. Then any element of $G$ can be written in the form $t^{-p}ht^{q}$ with $p,q\geq 0$ and $h \in H$. Moreover the normal closure $\langle H\rangle^G$ of $H$ in $G$ is given by $\langle H\rangle^G = \bigcup_{i\geq 0} t^{-i}Ht^i$.
\end{lem}

\Proof The claim about the form elements of $G$ take follows since for any $h \in H$, $th = \theta(h) t$ and similarly $h t^{-1} = t^{-1} \theta(h)$ in $G$. For the second part, notice that $t^{-i}Ht^i \subset \langle H\rangle^G$ for any $i$. Since $G/\langle H\rangle^G \cong \langle t\rangle$, we have that if $g = t^{-p}ht^q \in \langle H\rangle^G$, then $p=q$. Thus $\langle H\rangle^G = \bigcup_{i\geq 0} t^{-i}Ht^i$. \qed

\begin{lem} \label{productprimitive}
Let $G = \langle F,t \mid txt^{-1} =\theta(x), x\in F \rangle$ be any ascending HNN extension of a free group $F$ of finite rank with $\im(\theta) \subgroup [F,F]$. Suppose that $x, y \in F$ are non-primitive in $G$ and generate a free subgroup of rank $2$. Then $xy$ is primitive in $G$.
\end{lem}

\Proof
Suppose $x$, $y$ and $xy$ are all non-primitive. Let $x = u^m,y=v^n,xy= w^l$ for some $u,v,w \in G$ and $m,n,l\geq 2$. Let $q$ be the canonical quotient map from $G$ to $\langle t\rangle \cong \BZ$ mapping $F$ to $0$. Then $u,v,w$ lie in the kernel since $x$ and $y$ lie in the kernel. Note that the kernel is just the normalizer of $F$ in $G$. By \Cref{HNNextensionkernel}, there exist some $p \geq 0$ such that $u,v,w$ lie in the free subgroup $t^{-p}Ft^p$. But by \cite{LS62}, the equation $u^mv^n=w^l$ has a solution in a free group only if $u,v,w$ generate a free subgroup of rank $1$. This contradicts our hypothesis on $x$ and $y$.
\qed

We need the following lemma.

\begin{lem}\label{abelianprime}
Let $f \colon A \rightarrow A$ be an automorphism of a free abelian group $A$. If $f (ka) = la$ for some $a \neq 0$ and positive integers $k,l$, then $k=l$.
\end{lem}

\Proof We can assume without loss of generality that $A$ has infinite rank, so $A \cong \bigoplus_{i \in I} \BZ$ for some infinite index set $I$. We call a non-trivial element $a \in A$ prime if the common divisor of its finitely many non-zero coordinates is trivial. Note that any non-trivial $a \in A$ can be written as $a = d \cdot x$ with $x$ being prime. Since $f$ is an automorphism, it will preserve prime elements. So now suppose that $f(ka) = la$ with $k, l \in \BN$ and $a \neq 0$. We write $a = d x$ as above with $x$ prime. Then $k f(x) = l x$ and by cancelling common factors we might as well assume that $k$ and $l$ are coprime. Since $k$ divides all coordinates of the prime element $x$ it has to equal to one and the same holds for $l$ since $f(x)$ is prime.\qed

\begin{prop}\label{HNNcommutatorbvc}
Let $G = \langle F,t \mid txt^{-1} =\theta(x), x\in F \rangle$ be an ascending HNN extension of a free group $F$ of finite rank, where the image of $\theta$ lies in the commutator subgroup of $F$. If $x,y \in F \setminus [F,F]$ generate a free subgroup of rank $2$ in $F$ and $x$ is primitive, then the elements $\{x^kyx^ky^{-1} \mid k \geq 2 \}$ form pairwise distinct primitive conjugacy classes. In particular, $G$ does not have BVC.
\end{prop}

\Proof   Note first that $x^kyx^ky^{-1}$ does not lie in $[F,F]$ and is primitive when $k \geq 2$ by \Cref{productprimitive}. Note that every element in $G$ can be written in the form $t^{-p}wt^q$ for some $p \geq 0,q \geq 0$ and $w \in F_n$ by \Cref{HNNextensionkernel}. Now if $x^kyx^ky^{-1}$ is conjugate to $x^lyx^ly^{-1}$ for some $k \neq l$, then $x^kyx^ky^{-1} = t^{-p}wt^q  x^lyx^ly^{-1} t^{-q}w^{-1}t^p$ for some $p,q \geq 0$ and $w \in F$.  Hence $\theta^p(x^kyx^ky^{-1}) = w\theta^q(x^lyx^ly^{-1})w^{-1}$

If $p\neq q$, the equation never holds. In fact, assume $p>q$. We can further assume $\theta^q(x) \in \Gamma_r(F_n) \setminus \Gamma_{r+1}(F)$ for some $r \geq 2$ by \Cref{centralseries} (b). Then $\theta^p(x)\in \Gamma_{r+1}(F)$ by \Cref{derivedrank} and thus $\theta^p(x^kyx^ky^{-1}) \in \Gamma_{r+1}(F)$. On the other hand, $\theta^q(x^l) \in \Gamma_r(F) \setminus \Gamma_{r+1}(F)$ for any $l>0$ since $\Gamma_{r}(F)/\Gamma_{r+1}{F}$ is a free abelian group by \Cref{centralseries} (c). Now $\theta^q(x^lyx^ly^{-1}) = \theta^q(x^{2l}) [\theta^q(x^{-l}),\theta^q(y)] $ and $[\theta^q(x^{-l}),\theta^q(y)]  \in \Gamma_{r+1}(F)$ by \Cref{derivedrank}, we have $\theta^q(x^lyx^ly^{-1}) \in \Gamma_{r}(F) \setminus \Gamma_{r+1}(F)$. But the right hand side $w\theta^q(x^lyx^ly^{-1})w^{-1} \in  \Gamma_{r}(F)\setminus \Gamma_{r+1}(F)$, hence the equation cannot hold.

If $p=q$, then the equation again cannot hold unless $k=l$. In fact, assume $\theta^p(x) \in \Gamma_r(F) \setminus \Gamma_{r+1}(F)$, then both sides lie in $\Gamma_r(F) \setminus \Gamma_{r+1}(F)$ by the same argument above. 
By taking the quotient by $\Gamma_{r+1}(F)$ we obtain an equation in the free abelian group $\Gamma_r(F) / \Gamma_{r+1}(F)$. Then we would have $k([\theta^p(x)] + [\theta^p(yxy^{-1})]) = l([w\theta^p(x)w^{-1}] + [w\theta^p(yxy^{-1})w^{-1}])$. Note that $\Gamma_r(F) / \Gamma_{r+1}(F)$ is a free abelian group of infinite rank by \Cref{centralseries} (c) and the action of $w$ on $\Gamma_r(F) / \Gamma_{r+1}(F)$ induced by conjugation is an isomorphism. Thus the equation $k([\theta^p(x)] + [\theta^p(yxy^{-1})]) = l([w\theta^p(x)w^{-1}] + [w\theta^p(yxy^{-1})w^{-1}])$ can never hold unless $k=l$ by \Cref{abelianprime}. \qed

We are ready to prove the following.

\begin{thm} \label{HNN-ext-of-free-groups-BVC}
Let $G$ be an HNN extension of a free group of finite rank, then $G$ does not have BVC.
\end{thm}

\Proof By \Cref{HNNindexbvc}, we can assume  $G = \langle F_n,t \mid txt^{-1} =\theta(x), x\in F_n \rangle$, where $\theta \colon F_n \rightarrow F_n$ is injective, $F_n$ is a free group of rank $n$. For $n=1$ the group $G$ is solvable but not virtually cyclic. Thus $G$ does not have BVC by \Cref{BVCvirtuallysolvable}. So in the following we assume that $F_n$ is a free group of rank bigger than~$1$.

Note first that we have an induced map $\bar{\theta} \colon F_n/[F_n,F_n] \rightarrow F_n/[F_n,F_n] \cong \BZ^n$. Since the rank of the abelian group is finite, there exists some $k \geq 1$ such that $\rank( \ker(\bar{\theta}^{k+1})) = \rank( \ker(\bar{\theta}^k))$. But since $\BZ^n$ is free abelian, it follows that $\ker( \bar{\theta}^{k+1} ) = \ker(\bar{\theta}^{k})$, and we will denote this group by $K$. This implies that $\bar{\theta}^k$ induces an injective endomorphism of $\BZ^n/K$. If $K$ is a proper subgroup of $\BZ^n$, we consider the induced quotient map $F_n *_{\theta^k} \to (\BZ^n/K)*_{\bar{\theta}^k}$. Note that the quotient is a torsion-free metabelian group which is not virtually cyclic. Hence $F_n *_{\theta^k}$ does not have BVC by \Cref{BVCvirtuallysolvable} and \Cref{BVCquotient}. As $F_n *_{\theta^k}$ is a finite index subgroup of $F_n *_{\theta}$ (see for example \cite[2.2]{Ka}) we conclude that the latter group does not have BVC by \Cref{finiteindexsubgroupbvc}.

If $K = \BZ^n$, we are in the situation that the image of $\theta^k$ lies in the commutator subgroup of $F_n$. By \Cref{HNNcommutatorbvc} the group $F_n *_{\theta^k}$ does not have BVC. Again by \Cref{finiteindexsubgroupbvc} it follows that $F_n *_{\theta}$ does not have BVC.
\qed

We now want to apply the previous results to verify Conjecture B for the class of one-relator groups. Recall that a \textbf{one-relator group} is a group $G$ which has a presentation with a single relation, so $G = \langle x_1, \ldots, x_n \mid r \rangle$ where $r$ is a word in the free group $F$ on the letters $x_1, \ldots, x_n$. The group $G$ is torsion-free precisely when $r$, as an element of the free group $F$, is not a proper power. If $r = s^n$ for some maximal $n \geq 2$ and $s \in F$, then $s$, considered as an element in $G$, is of order $n$. In all cases there exists a finite $G$-CW model for $\u E G$, see for example \cite[4.12]{Lu05}.

For one-relator groups $G$ with torsion, Newman's Spelling Theorem \cite{Ne} implies that $G$ is a hyperbolic group. In particular, the one-relator groups containing torsion satisfy Conjecture B. However, our proof of the following theorem does not depend on this fact.

\begin{thm}\label{one-relator-BVC}
A one-relator group has BVC if and only if it is cyclic.
\end{thm}

\Proof Let $G$ be a one-relator group.

If the one-relator presentation of $G$ contained three or more generators then $G$ would surject to $\BZ^2$, in particular $G$ would not have BVC by \Cref{BVChomology}.
Thus we can restrict to the case that $G$ has two generators, so
$$
G = \langle a, b \mid R(a,b) = 1 \rangle
$$
for some word $R(a,b)$ in the free group on the two generators $a,b$. By \cite[Lemma V.11.8]{LS77} we can moreover assume that the exponent sum of one of the generators in the single relator equals to zero, say for the generator $a$. The following rewriting procedure, which we outline for the reader's convenience, is standard. The proofs of the mentioned facts can be found in \cite[IV.5]{LS77}. We let $b_i = a^{i} b a^{-i}$ for all $i \in \BZ$. Then $R$ can be rewritten as a cyclically reduced word $R'$ in terms of these, so $R' = R'(b_m, \ldots, b_M)$ for some $m \leq M$, such that the elements $b_m$, $b_M$ occur in $R'$. If $m = M$, then $R(a,b) = b^m$ for some $m \in \BZ$ and thus $G \cong \BZ$ or $G \cong \BZ * \BZ/|m|$ where $|m| \geq 2$. Note that by \Cref{free-products-bvc} the latter group does not have BVC. So in the following we can assume that $m < M$. We let
$$
H = \langle b_m, \ldots, b_M \mid R'(b_m, \ldots, b_M) = 1 \rangle.
$$
Moreover we define $A$ to be the subgroup of $H$ generated by $b_m, \ldots b_{M-1}$ and we let $B$ to be the subgroup of $H$ generated by $b_{m+1}, \ldots b_M$. Then $A$ and $B$ are free subgroups of the one-relator group $H$ and $G$ is isomorphic to the HNN extension $H*_{\theta}$ where $\theta \colon A \to B$ is the isomorphism defined by $\theta(b_i) = b_{i+1}$ for $m \leq i < M$.

If $[H:A] \geq 2$ and $[H:B] \geq 2$, then $G$ does not have BVC by \Cref{HNNindexbvc}. Otherwise $G$ is an ascending HNN extension, say with $H = A$. Since $A$ was free, $G$ is an ascending HNN extension of a finitely generated free group. The claim now follows from \Cref{HNN-ext-of-free-groups-BVC}.
\qed

\section{Groups with some hyperbolicity} \label{hyper}
In this section, we first show that acylindrically hyperbolic groups do not have BVC. Using this we show that any finitely generated 3-manifold group does not have BVC.

We first give a quick definition of acylindrically hyperbolic group and refer to \cite{Hu} and \cite{Os16}   for more information. Recall the action of a group $G$ on a metric space $S$ is called \textbf{acylindrical} if for every $\varepsilon >0 $ there exist
$R, N >0$ such that for every two points $x, y$ with $d(x, y) \geq R$, there are at most $N$ elements $g \in G$ satisfying $ d(x,gx) \leq \varepsilon $ ~and~ $d(y,gy) \leq \varepsilon $. Given a hyperbolic space $S$, we use $\partial S$ to denote its Gromov boundary.

\begin{defn}
A group $G$ is called \textbf{acylindrically hyperbolic} if there exists a generating set $X$ of $G$ such that the corresponding Cayley graph $\Gamma(G, X)$
is hyperbolic,  $|\partial \Gamma(G, X)| > 2$, and the natural action of G on $\Gamma(G, X)$ is acylindrical.
\end{defn}

\begin{prop} \label{acylindrically hyperbolic}
An acylindrically hyperbolic group does not have BVC.
\end{prop}

\Proof Let $G$ be an acylindrically hyperbolic group. By \cite[Theorem 1.2]{Os16} this is equivalent to saying that $G$ contains a proper infinite hyperbolically embedded subgroup. By \cite[Lemma 6.14]{DGO}, we further have a subgroup $K = F_2 \times K(G)$ inside $G$ which is hyperbolically embedded, where $K(G)$ is some maximal normal finite subgroup of $G$ and $F_2$ is the free group of rank $2$. The following two statements are copied from the proof of \cite[VI.1.1]{Hu} (they can be easily decduced from \cite[Proposition 4.33]{DGO}): (1) every element of infinite order which is primitive in $K$ is also primitive in $G$; (2), if two elements of $K$ of infinite order are conjugate in $G$, they are conjugate in $K$. Now since $K$ contains $F_2$ as a direct factor, there exists infinitely many primitive elements of infinite order $g_1,g_2,\ldots,g_n,\ldots$ in $F_2 \subset K$ such that $g_i$ is not conjugate to $g_j$ or $g_j^{-1}$ when $i \neq j$.  But any two such elements or conjugates of them cannot lie in the same virtually cyclic subgroup. In fact, let $g_i,g_j$ be two such elements and suppose $g_i$, $xg_jx^{-1}$ lie in the same virtually cyclic subgroup in $G$, where $x \in G$. Then we have $g_i^m = xg_j^nx^{-1}$ for some $m$ and $n$ by \Cref{vircycorder}. Thus $g_i^m$ and $g_j^n$ are conjugate in $G$. Hence they are also conjugate in $F_2$, which is a contradiction. Thus these primitive elements of infinite order cannot lie in finitely many conjugacy classes of virtually cyclic subgroups. Hence $G$ does not have BVC.

\qed

\begin{cor}\label{acc-exp}
The following classes of groups do not have BVC:
\begin{enumerate}
    \item Hyperbolic groups that are not virtually cyclic;
    \item Non-elementary groups that are hyperbolic relative to proper subgroups;
    \item The mapping class group $MCG(\Sigma_g,p)$ of a closed surface of genus $g$ with $p$ punctures except for $g = 0$ and $p \leq 3$ (in these exceptional cases, $MCG(\Sigma_g,p)$ is finite);
    \item $Out(F_n)$ where $n \geq 2$;
    \item Groups which act properly on proper CAT(0) spaces and contain rank-one elements;
\end{enumerate}
\end{cor}
\Proof These groups are all acylindrically hyperbolic, we refer to \cite[I.1,page 4]{Hu} and \cite[Section 8]{Os16} for the detailed reference. \qed

\begin{cor} \label{surjecachyp}
Let $\phi \colon A \rightarrow G$ be a surjective group homomorphism and suppose that $G$ is an acylindrically hyperbolic group. Then $A$ does not have BVC.
\end{cor}

\Proof As in the proof of \Cref{acylindrically hyperbolic}, there exists infinitely many primitive conjugacy classes of elements $\{g_i\}_{i\geq 1}$ inside a certain free group $F_2 \subset G$ such that ${g_i}^{p_i}$ is not conjugate to any ${g_j}^{p_j}$ for $i\neq j$ and $p_i,p_j \neq 0$. Now we take any preimage  $g_i'$ of $g_i$ in $A$. Then the $g_i'$ are primitive and ${g'_i}^{p'_i}$ is not conjugate to any ${g'_j}^{p'_j}$ for $i\neq j$ and $p'_i,p'_j \neq 0$. By \Cref{vircycorder}, they cannot lie in finitely many conjugacy classes of virtually cyclic subgroups in $A$. Hence $A$ does not have BVC. \qed

By a \textbf{$3$-manifold group} we mean a group that can be realized as the fundamental group of a $3$-manifold, which may be open or have boundary. Note that, by Scott's theorem \cite{Sc}, every finitely generated $3$-manifold group is itself the fundamental group of a compact 3-manifold. Minsanyan and Osin prove in \cite[2.8]{MO} the following:
\begin{lem} Let $G$ be a subgroup of the fundamental group of a compact $3$-manifold, then exactly one of the following holds:
\begin{enumerate}
    \item $G$ is acylindrically hyperbolic.
    \item $G$ contains a normal infinite cyclic subgroup $N$ such that $G/N$ is acylindrically hyperbolic.
    \item $G$ is virtually polycyclic.
\end{enumerate}
\end{lem}

Now combining this with \Cref{BVCvirtuallysolvable}, \Cref{acylindrically hyperbolic} and \Cref{surjecachyp}, we have the following
\begin{prop}\label{3-manifold-group-BVC}
Let $G$ be the subgroup of the fundamental group of a compact $3$-manifold, then $G$ has BVC if and only if $G$ is virtually cyclic. In particular, if $G$ is a finitely generated $3$-manifold group, then $G$ has BVC if and only if $G$ is virtually cyclic.
\end{prop}

\section{Groups acting on CAT(0) spaces}\label{CAT(0)}

In this section, we study groups acting on CAT(0) spaces and show that many of them do not have BVC.  We first give a quick review of properties of CAT(0) spaces and groups that we may need and refer to \cite{BH} for more details.

\begin{defns}\cite[II.6.1]{BH} Let $X$ be a metric space and let $g$ be an isometry of $X$. The \textbf{displacement function} of $g$ is the function $d_g \colon X \rightarrow \BR_+ =\{r\geq 0 \mid r\in \BR\}$ defined by $d_g(x) =d(g x,x)$. The \textbf{translation length} of $g$ is the number $|g|:= \inf \{d_g (x) \mid x\in X\}$. The set of points where $d_g$ attains this infimum will be denoted by $Min(g)$. More generally, if $G$ is a group acting by isometries on $X$, then $Min(G):=\bigcap_{g\in G} Min(g)$. An isometry~$g$ is called \textbf{semi-simple} if $Min(g)$ is non-empty. An action of a group by isometries of $X$ is called \textbf{semi-simple} if all of  its elements are semi-simple.
\end{defns}

The following theorem is known as the Flat Torus Theorem \cite[II.7.1]{BH}.
\begin{thm} \label{flattorusthm}
Let $A$ be a free abelian group of rank $n$ acting properly by semi-simple isometries on a CAT(0) space $X$. Then:

\begin{enumerate}
\item $Min(A) = \bigcap_{\alpha\in A} Min(\alpha)$ is non-empty and splits as a product $Y \times \BE^n$, here $\BE^n$ denotes $\BR^n$ equipped with the standard Euclidean metric.
\item Every element $\alpha \in A$ leaves $Min(A)$ invariant and respects the product decomposition; $\alpha$ acts as the identity on the first factor $Y$ and as a translation on the second factor $\BE^n$.
\item The quotient of each $n$-flat $\{y\} \times \BE^n$ by the action of $A$ is an $n$-torus.
\end{enumerate}
\end{thm}

It is clear that the translation length is invariant under conjugation, i.e. $|h g h^{-1}| = |g|$ for any $g,h \in G$. Moreover, for $g$ semi-simple, we have that $|g^n| = |n| \cdot |g|$ for any $n \in \BZ$, e.g. by the Flat Torus Theorem. It turns out that the translation length can also be defined by the following limit for $g$ a semi-simple isometry
$$
|g| = \lim_{n \to \infty} \frac{1}{n} d(x, g^n x),
$$
where $x$ is an arbitrary point of the CAT(0) space $X$ \cite[II.6.6]{BH}.

Note that if a group acts properly and cocompactly on a CAT(0) space via isometries, we call the group a \textbf{CAT(0) group}. In this case, the action is semi-simple.

\begin{prop}  \label{prop:CAT(0)}
If a group $G$ acts properly and cocompactly by isometries on a CAT(0) space $X$, then
\begin{enumerate}

    \item $G$ is finitely presented.
    \item $G$ has only finitely many conjugacy classes of finite subgroups.
    \item Every solvable subgroup of $G$ is virtually abelian.
    \item Virtually abelian subgroups of $G$ satisfy the ascending chain condition.

    \item $G$ has a finite model for $\u E G$.

    \item  There is a finite-dimensional model for $\uu E G$.

\end{enumerate}

\end{prop}

\Proof (a) - (c) can be found in \cite[III.$\Gamma$.1.1]{BH}, (d) can be found in \cite[II.7.5]{BH}. Since $G$ acts on $X$ properly and cocompactly via isometry, $X$ is proper by \cite[I.8.4(1)]{BH}. With this (e) is implied by  \cite[Proposition A]{On}. The last statement was proven in \cite{Lu09}. Also Farley has given a construction of $\uu EG$ in \cite{DF10} for some CAT(0) groups, however without controlling the dimension. \qed


\begin{lem}\label{vctranslationlength}
Let $V$ be an infinite virtually cyclic group which acts on a CAT(0) space via semi-simple isometries. Then there exist an element $v_0 \in V$ such that for any element $v \in V$, the translation length  $|v|$ of $v$ is an integer multiple of the translation length of $v_0$.
\end{lem}

\Proof When $v$ is of finite order, then $|v| = 0$. So let us assume in the following that $v$ is of infinite order, in this case the translation length is strictly positive.
If $v^k = w^{k'}$ for some $w \in V$ and $k, k' \in \BZ$, then $|k| |v| = |v^k| = |w^{k'}| = |k'| |w|$. Now by \Cref{vircycorder}, there exist $v_0 \in V$ such that for any $v \in V$, there exist nonzero $p,p_0$ such that $v_0^{p_0} = v^p$ with $\frac{p_0}{p} \in \BZ$. This implies that $|v|$ is a multiple of $|v_0|$.

\qed

The lemma leads us to define the following terminology.

\begin{defn}
We define a subset $A$ of the real numbers to be \textbf{finitely divisor dominated} if there are finitely many real numbers $x_1,x_2,\ldots, x_n$ such that every $a \in A$ can be written in the form $kx_i$ for some $k\in \BZ$, or equivalently
$$
A \subset \bigcup_{i = 1}^n \BZ \cdot x_i.
$$

In this case we say that $A$ is finitely divisor dominated by $x_1,x_2,\ldots, x_n$.
\end{defn}

Note that for a CAT(0) group the set of translation lengths $\{ \, |g| \mid g \in G \}$ is a discrete subset of $\BR$ \cite[II.6.10 (3)]{BH}. We obtain first the key property of the set of translation lengths for a group acting on a CAT(0) space with BVC.

\begin{lem}\label{BVCfdd}
Let $G$ be a group acting properly on a CAT(0) space via semi-simple isometries.
Let $L = \{ \, |g| \mid g \in G \} \subset \BR_{\geq 0}$ be the set of translation lengths of $G$. If $G$ has BVC, then $L$ is finitely divisor dominated.
\end{lem}

\Proof Note first that if $g$ and $h$ are conjugate, then they have the same translation length. Now assume $G$ has BVC, and let $V_1,V_2,\ldots,V_n$ be witnesses. We only need to consider those $V_i$ that are infinite since a finite order element has vanishing translation length. By \Cref{vctranslationlength}, we can choose for each $V_i$ an element $v_i$ such that the translation length of any infinite order element of $V_i$ is a multiple of the translation length $|v_i|$. Now $L$ is finitely divisor dominated by $|v_1|,|v_2|,\ldots, |v_n|$.
\qed

\begin{rem}
For a hyperbolic group acting on its Cayley graph with respect to some fixed generating set $S$, we define the algebraic translation length using the limit  $
|g| = \lim_{n \to \infty} \frac{1}{n} d_S(1, g^n )
$ \cite[III.$\Gamma$.3.13] {BH}, where $d_S$ denotes the word metric with respect to $S$. Gromov \cite[III.$\Gamma$.3.17]{BH} showed that the set of algebraic translation lengths in this case is a discrete subset of the rational numbers and the denominators are bounded. In particular, the set of algebraic translation lengths of a hyperbolic group is finitely divisor dominated.
\end{rem}

We need do some preparation before we prove our main result in this section.

\begin{lem} \label{fddrational1}
Let $x >0,y\geq 0$ be two rational numbers and let
$$A = \left \{\lambda_n \mid \lambda_n = \sqrt {x + (y+n)^2 }  ,n \in \BN \right \}.$$
Then $A$ is not finitely divisor dominated.
\end{lem}

\Proof Since $x$ and $y$ are rational numbers, we can choose $d$ to be the smallest positive integer such that $2yd$ and $(y^2+x)d$ are integers. Then we can consider the quadratic polynomial
$$f(n) =  d(x+(y+n)^2) = dn^2 + 2ydn +(y^2 + x)d$$
which has coprime integer coefficients. Note that $f(n)$ is irreducible over $\BR$ since $d,x>0$. Now by an old result of Ricci \cite{Ri}, there exists infinitely many positive integers $n$ such that the integer $f(n)$ is square-free.

Now if $A$ was finitely divisor dominated, there would exist finitely many positive real numbers $x_1,x_2,\ldots,x_m$ such
$$
A \subset \bigcup_{i = 1}^n \BZ \cdot x_i
$$

In particular,  there exists some $i_0$ and an infinite sequence $n_1,n_2,\ldots,n_j,\ldots $ of natural numbers such that $\lambda_{n_j} = k_j x_{i_0}$, with $k_j \in \BZ $ and $f(n_j)$ square-free. This implies there are infinitely many $k_j > k_1$ such that
$$\frac{\lambda_{n_j}^2}{\lambda_{n_1}^2} = \frac{x + (y+n_j)^2}{x + (y+n_1)^2} = \frac{k_j^2}{k_1^2}$$
This further implies
$$ n_j^2 + 2yn_j +y^2 + x = \frac{k_j^2}{k_1^2} (n_1^2 + 2yn_1 +y^2 + x) $$
Multiplying both sides by $d$, we obtain
$$dn_j^2 + 2ydn_j +(y^2 + x)d = \frac{k_j^2}{k_1^2} (dn_1^2 + 2ydn_1 +(y^2 + x)d) $$
Now since $f(n) = d(x+(y+n)^2) = dn^2 + 2ydn +(y^2 + x)d$
is a polynomial in $n$ with positive integer coefficients that have no common divisor, the left hand side of the above equation must be a positive integer and $f(n_1) = dn_1^2 + 2ydn_1 +(y^2 + x)d$ is also a positive integer. But since $k_j > k_1$ are positve integers, the value of $f(n_j) = \frac{k_j^2}{k_1^2} (dn_1^2 + 2ydn_1 +(y^2 + x)d)$ is not square-free. This leads to a contradiction as we have chosen the $n_j$ such that $f(n_j)$ is square-free.
\qed

\begin{lem} \label{fddrational2}
Let $n_0,q_0 \in \BN$ and $x>0, y\geq 0$ be two real numbers, such that there are infinitely many integers $m>n_0$ with $\frac{x^2 + (y+m)^2}{x^2 + (y+n_0)^2} = \frac{p^2}{q_0^2}$ for some $p \in \BN$. Then $y \in \BQ$, $x^2 \in \BQ$.

\end{lem}

\Proof Let $m_1,m_2,\ldots, m_i, \ldots$, be an infinite sequence of positive integers greater than $n_0$ such that $\frac{x^2 + (y+m_i)^2}{x^2 + (y+n_0)^2} = \frac{p_i^2}{q_0^2}$. Let us write $x^2 + (y+n_0)^2 = q_0^2t$ for some $t > 0$, then $x^2 + (y+m_i)^2 = p_i^2t$. Subtracting this by the previous equation, we get
$$ (m_i - n_0)(2y+m_i +n_0) = (p_i^2-q_0^2)t.$$
Now comparing this with the same equality for $m_j$, we obtain
\begin{equation} \label{equation}
\frac{ (m_i - n_0)(2y+m_i +n_0)}{(m_j - n_0)(2y+m_j +n_0)} = \frac{p_i^2-q_0^2}{p_j^2-q_0^2}
\end{equation}
 Since $m_i,m_j,n_0,p_i,p_j$ are all integers, we have $y$ is rational unless

 $$
\frac{p_i^2-q_0^2}{p_j^2 - q_0^2} = \frac{m_i^2-n_0^2}{m_j^2-n_0^2}
$$

But this cannot happen as long as $m_i \neq m_j$. In fact, note first that we can assume without loss of generality that $n_0 = 0$ (via shifting $y$ by some integer). Now let $r = (p_i^2-q_0^2)/(p_j^2-q_0^2)$. Then equation \ref{equation} above leads to
$$
2y(m_i - m_j r) = r m_j^2 - m_i^2
$$
We cannot solve for $y$ if $m_i = m_j r$. But if this happens $r = \frac{m_i^2}{m_j^2} = \frac{m_i}{m_j}$,  hence $m_i = m_j$.

 This also immediately implies $x^2 \in \BQ$ using the equation $x^2 + (y+n_0)^2 = q_0^2t$.
\qed

Combining \Cref{fddrational1} and \Cref{fddrational2}, we have the following

\begin{cor}\label{lengthFM}
Let $x>0,y  \geq 0$ be two real numbers and let $A = \{\lambda_n \mid \lambda_n = \sqrt {x^2 + (y+n)^2 }  ,n \in \BN \}$. Then $A$ is not finitely divisor dominated.
\end{cor}

\begin{prop} \label{abeliansubgroup}
Let $G$ be a group acting properly on a CAT(0) space $X$ via semi-simple isometries. If $G$ contains $\BZ^2$ as a subgroup, it does not have BVC.
\end{prop}
\Proof Assume $G$ has BVC, then the set of translation lengths $L = \{ \,|g| \mid g \in G \}$ is finitely divisor dominated by \Cref{BVCfdd}. On the other hand, by the Flat Torus Theorem, we have that $\BZ^2$ acts on a flat plane $\BE^2$ inside $X$ and the translation length of any $g = (z,w) \in \BZ^2$ is just $d(g,gx_0)$ for some base point $x_0 \in \BE^2$. Let $a$ be the translation vector of $(1,0) \in \BZ^2$ and $b$ for $(0,1) \in \BZ^2$. Let $(a_1,a_2), (b_1,b_2)$ be the coordinate of $a$ and $b$ in the Euclidean plane $\BE^2$. Without loss of generality, we can assume $a_1 >0 , a_2 \geq 0$ and $(b_1,b_2) = (1,0)$. Then the translaton length of $(1,k) \in \BZ^2$ is the length of the vector $a +kb =(a_1,a_2) + k(0,1) = (a_1,a_2+k)$, which is

$$\lambda_{(1,k)} = \sqrt{ a_1^2+(a_2+k)^2}.$$

Now if a set is finitely divisor dominated, then any subset of it is also finitely divisor dominated. In particular, the subset $\left\{\sqrt{ a_1^2+(a_2+k)^2} \mid k \in \BN \right \}$ is finitely divisor dominated for some $a_1 >0,a_2\geq 0$. But this contradicts \Cref{lengthFM}.

\qed

Recall that a semi-simple isometry is called \textbf{hyperbolic} if it has positive translation length. Now if $g$ acts properly on a CAT(0) space $X$ via a hyperbolic isometry, by the Flat Torus Theorem, we have an axis $\BE^1$ on which $g$ acts via translation by the length $|g|$.

\begin{defn}
Supppose $g$ is a hyperbolic isometry of a CAT(0) space $X$. Then $g$ is called \textbf{rank one} if no axis of $g$ bounds a flat half plane in $X$.
\end{defn}

Note that if a group acts on a CAT(0) space $X$ properly and cocompactly via isometries, then $X$ is proper \cite[I.8.4(1)]{BH}. Combining this with \Cref{acc-exp} (e) and \Cref{abeliansubgroup}, we have the following

\begin{thm}\label{BVCCAT(0)}
Let $G$ be a subgroup of a CAT(0) group which contains a subgroup ismorphic to $\BZ^2$ or a rank-one isometry, then $G$ does not have BVC.
\end{thm}

\begin{rem}\label{flat-closing-conj}
The Flat Closing Conjecture \cite[6.B3]{Gr} predicts that $X$ contains
a $d$-dimensional flat if and only if $G$ contains a copy of $\BZ^d$ as a subgroup. In particular, it implies that a CAT(0) group is hyperbolic if and only if it does not contain a subgroup isomorphic to $\BZ^2$. Thus the Flat Closing Conjecture together with \Cref{BVCCAT(0)} would also imply that a CAT(0) group has BVC if and only if it is virtually cyclic.
\end{rem}

Recall that a \textbf{CAT(0) cube group} is a group which acts properly and cocompactly on a CAT(0) cube complex via isometries.

\begin{lem} \label{CAT(0)cuberankrigidity}
Let $G$ be a group which acts on a CAT(0) cube complex $X$ properly and cocompactly via isometries and suppose that $G$ is not virtually cyclic. Then either $G$ contains a rank one isometry or $G$ contains a free abelian subgroup of rank $2$.
\end{lem}
\Proof This is essentially due to Caprace and Sageev \cite{CS}.  Note first that $X$ is locally finite, see for example \cite[I.8.4(1)]{BH}. Note also that $G$ acts on $X$ without fixed points and essentially (see \cite[1.1]{CS} for the terminology). Now by  \cite[Theorem A]{CS} and remarks below it, we have that either
$X$ is a product of two unbounded CAT(0) cube subcomplexes or $G$ contains an element acting on $X$ as a rank one isometry. Note that since $G$ acts on $X$ cocompactly, if $X$ is a product of two CAT(0) cube complexes, by \cite[Corollary D]{CS}, it follows that $X$ contains a free abelian subgroup of rank $2$.
\qed

\begin{cor}\label{BVCCAT(0)cube}
Let $G$ be a CAT(0) cube group. Then $G$ has $BVC$ if and only if $G$ is virtually cyclic.
\end{cor}
\qed

\bibliographystyle{amsplain}

\end{document}